\documentclass[11pt]{article}

\usepackage{amssymb}
\usepackage{latexsym}
\usepackage[dvips]{epsfig}

\newcommand{\ben}{\begin{enumerate}}
\newcommand{\een}{\end{enumerate}}
\newcommand{\bde}{\begin{defn}}
\newcommand{\ede}{\end{defn}}
\newcommand{\bex}{\begin{exa}}
\newcommand{\eex}{\end{exa}}
\newcommand{\barr}{\begin{array}}
\newcommand{\earr}{\end{array}}
\newcommand{\btab}{\begin{tabular}}
\newcommand{\etab}{\end{tabular}}
\newcommand{\beq}{\begin{equation}}
\newcommand{\eeq}{\end{equation}}
\newcommand{\bea}{\begin{eqnarray*}}
\newcommand{\eea}{\end{eqnarray*}}
\newcommand{\bce}{\begin{center}}
\newcommand{\ece}{\end{center}}

\newcommand{\iso}{\cong}

\newcommand{\fl}[1]{\lfloor #1 \rfloor}

\newcommand{\flf}[2]{\left\lfloor\frac{#1}{#2}\right\rfloor}

\newcommand{\ree}[1]{(\ref{#1})}

\newcommand{\ra}{\rightarrow}

\newcommand{\qed}{\mbox{$\Diamond$}\vspace{\baselineskip}}

\newtheorem{theorem}{Theorem}[section]
\newtheorem{proposition}[theorem]{Proposition}
\newtheorem{lemma}[theorem]{Lemma}
\newtheorem{definition}[theorem]{Definition}
\newtheorem{corollary}[theorem]{Corollary}

\newtheorem{conjecture}[theorem]{Conjecture}
\newenvironment{proof}{\noindent {\bf Proof:}}{{\qed}}

\newenvironment{proofofstrongstrom}{\noindent {\bf Proof of Theorem~\ref
{strongstrom}:}}{{\qed}}

\newtheorem{thm}{Theorem}[section]

\newtheorem{exa}[thm]{Example}
\newtheorem{defn}[thm]{Definition}

\newcommand{\ble}{\begin{lemma}}
\newcommand{\ele}{\end{lemma}}
\newcommand{\bth}{\begin{thm}}
\renewcommand{\eth}{\end{thm}}
\newcommand{\bpr}{\begin{proposition}}
\newcommand{\epr}{\end{proposition}}
\newcommand{\bco}{\begin{corollary}}
\newcommand{\eco}{\end{corollary}}

\newcommand{\vanish}[1]{}
\newcommand*{\newnot}{\mathop{\mathrm{not }}}
\newcommand*{\achain}{\mathcal{A}}
\newcommand*{\minusmax}{\overline}
\newcommand*{\pat}{\achain_1\oplus\achain_l}
\newcommand*{\packpat}{c_{Q_l}}
\newcommand*{\packpatabbrev}{c}

\newcommand*{\pxy}{P^{x\rightarrow y}}
\newcommand*{\permpat}{1(l+1)l\dots 2}

\begin{document}
\title{Pattern frequency sequences and internal zeros}

\author{Mikl\'os B\'ona
\thanks{University of Florida, Gainesville FL 32611.
Email: bona@math.ufl.edu}
\and Bruce E. Sagan \thanks{Michigan State University, E. Lansing, MI 48824.
Email: sagan@math.msu.edu}
\and Vincent R. Vatter \thanks{Michigan State University,  E. Lansing, MI
48824. Email: vattervi@msu.edu }\\[5pt]
This paper is dedicated to the memory of Rodica Simion\\
who did some seminal work in the area of pattern avoidance}

 \date{}

\maketitle

\begin{abstract}
Let $q$ be a pattern and let $S_{n,q}(c)$ be the number of
$n$-permutations having exactly $c$ copies of $q$.  We investigate
when the sequence $(S_{n,q}(c))_{c\ge0}$ has internal zeros.
If $q$ is a monotone pattern it turns out that, except for $q=12$ or $21$,
the nontrivial sequences (those where $n$ is at least the length of
$q$) always have internal zeros.
For the pattern $q=1(l+1)l\ldots2$
there are infinitely many sequences which contain internal zeros
and when $l=2$ there are also
infinitely many which do not.  In the latter case, the only possible
places for internal zeros are the next-to-last or the second-to-last
positions.  Note that by symmetry this completely determines the existence
of
internal zeros for all patterns of length at most three.
\end{abstract}

\section{Introduction}
Let $q=q_1 q_2\ldots q_l$ be a permutation in the
symmetric group $S_l$.  We call $l$ the {\it length} of $q$.
We say that the permutation
$p=p_1 p_2 \ldots p_n\in
 S_n$ {\it contains a $q$-pattern\/}
 if and only if there is a subsequence
$p_{i_1}p_{i_2}\ldots p_{i_l}$  of $p$ whose elements are in the same
relative order as those in $q$, i.e.,
$$
p_{i_j}<p_{i_k} \mbox{ if and only if } q_j<q_k
$$
whenever $1\le j,k \le l$.
For example, 41523 contains exactly two  132-patterns, namely 152 and
153.  We let
$$
\mbox{$c_q(p) =$ the number of copies of $q$ in $p$},
$$
so that
$c_{132}(41523)=2$.
Permutations  containing a given number of $q$-patterns have been
extensively studied
recently~[1--11].

In this paper, we consider permutations with a given number of $q$-patterns
from a new angle.  Let
$$
\mbox{$S_{n,q}(c) =$ the number of $n$-permutations with exactly $c$
patterns 
of type $q$}.
$$
For $n$ and $q$ fixed, the sequence $(S_{n,q}(c))_{c\ge0}$ is called the
{\em 
frequency sequence} of the pattern $q$ for $n$.  Clearly this sequence
consists entirely of zeros if $n$ is less than the length of $q$ and
so we call these sequences {\em trivial\/} and all others {\em nontrivial}.
We also say that an $n$-permutation $p$ is {\it $q$-optimal} if there
is no $n$-permutation with more  copies of $q$ than $p$, and let
$$
\mbox{$M_{n,q}=c_q(p)$ for an optimal $p$}.
$$

The only $q$ for which the frequency sequence is well understood is $q=21$
(or equivalently $q=12$). Occurences of this pattern are called inversions.
It is well known \cite{stanley} that for all $n$, the frequency sequence of
inversions is log-concave, and so is unimodal and has no internal zeros.

When $q$ is has length greater than 2, numerical evidence suggests
that the frequency
sequence of $q$ will no longer be unimodal, let alone log-concave. In fact,
internal zeros seem to be present in
most frequency sequences. An integer $c$ is called an {\it internal
zero\/} of the
sequence $(S_{n,q}(c))_{c\ge0}$ if for some $c$ we have
$S_{n,q}(c)=0$, but there exist $c_1$ and $c_2$
with  $c_1<c<c_2$ and $S_{n,q}(c_1)$, $S_{n,q}(c_2)\neq0$.

In the rest of this paper we study the frequency sequences of the
monotone pattern $q=12\ldots l$ and the pattern $q=1(l+1)l\ldots 2$.  We
will
show that in the first case, when $l\ge3$ (the case $l=2$ has already
been mentioned) the nontrivial sequences always have internal zeros.
For $1(l+1)l\ldots 2$-patterns there are infinitely many $n$ where the
sequence has internal zeros.  For the $132$-pattern
there are also infinitely many $n$ where the sequence has no internal zeros.
And internal
zeros can only appear in positions $M_{n,132}-1$ or $M_{n,132}-2$.


\section{The monotone case}

We will now consider the sequence $(S_{n,q}(c))_{c\ge0}$ where
$q=12\ldots l$.  For later reference, we single out the known case
when $l=2$ discussed in the introduction.
\begin{proposition}            \label{12}
The sequence $(S_{n,12}(c))_{c\ge0}$ has no internal zeros (and is, in
fact, log concave). The unique optimal permutation is $p=12\ldots n$
with
$$
M_{n,12}={n\choose 2}\quad\qed
$$
\end{proposition}

It turns out that this is the only monotone pattern (aside from 21)
whose sequence has no internal zeros.  To prove this result, define an
{\it inversion\/} (respectively, {\it noninversion}) in $p=p_1
p_2\ldots p_n$ to be a pair $(p_i,p_j)$
such that $i<j$ and $p_i>p_j$ (respectively, $p_i<p_j$).

\begin{theorem}                \label{12...l}
Let $q=12\ldots l$ where $l\ge3$.  Then in $S_n$, the unique optimal
permutation is $p=12\ldots n$ and
$$
M_{n,12\ldots l}={n\choose l}.
$$
The set of permutations having the
next greatest number of copies of $q$ are those obtained from $p$ by
an adjacent transposition and this number of copies is
\beq                                            \label{mono}
{n-1 \choose l}+{n-2\choose l-1}.
\eeq
\end{theorem}
\begin{proof}
Consider any $r\in S_n$ different from
$p$.  Then $r$ has an inversion $(r_i,r_j)$.  So the number of copies
of $q$ in $r$ is the number not containing $r_i$ plus the number which
do contain $r_i$.  The permutations in the latter case can not contain
$r_j$.  So~\ree{mono} gives an upper bound for the number of copies of
$q$ which is strict unless $r$ has exactly one inversion.  The theorem
follows.
\end{proof}

\bco
Let $q=12\ldots l$ where $l\ge3$.  Then for $n\ge l$ the sequence
$(S_{n,12\ldots l}(c))_{c\ge0}$ has internal zeros.
\eco
\begin{proof}
{From} the previous theorem, we see that the number of zeros directly
before $S_{n,q}(M_{n,q})=1$ is
$$
{n\choose l}-{n-1\choose l}-{n-2\choose l-1}={n-2\choose l-2}\ge n-2\ge1
$$
since $n\ge l\ge3$.
\end{proof}

For use in the 132 case, we record the following observation.
\ble                                                    \label{12/132}
For any integer $c$ with $0\le c\le {n\choose 2}$ there is a
permutation $p\in S_n$ having $c$ copies of the pattern $21$ and no copies
of
$132$.
\ele
\begin{proof}
We induct on $n$.  The result is clearly true if $n\le2$.  Assuming it
is true for $n-1$, first consider $c\le{n-1\choose 2}$ and let
$p\in S_{n-1}$ satisfy the lemma.  Then the concatenation $pn\in S_n$
works for such $c$.  On the other hand, if ${n-1\choose 2}<c\le {n\choose2}$
then consider $c'=c-(n-1)\le{n-1\choose2}$.  Pick $p\in S_{n-1}$ with
$c'$ copies of $21$ and none of $132$.  Then $np\in S_n$ is the
desired permutation.
\end{proof}

\section{The case $q=1(l+1)l\ldots2$ and layered patterns}

The rest of this paper is devoted to the study of the frequency sequences of
the patterns $\permpat$ for $l\ge 2$. To simplify notation,
and write $F_{n,\permpat}$ for the sequence
$(S_{n,\permpat}(c))_{c\ge0}$. One crucial property of these patterns
is that  
they are layered. This section gives an overview of some important results
on 
layered patterns.

A pattern is {\it layered\/} if it is the concatenation of subwords
(the {\em layers}) where the entries
decrease within each layer, and increase between the layers.
For example, $3\: 2\: 1\: 5\: 4\: 8\: 7\: 6\: 9$ is a layered pattern with
 layers $3\: 2\: 1,\: 5\: 4,\: 8\: 7\: 6$, and 9. Layered patterns are
examined in Stromquist's work~\cite{str:plp} and in Price's
thesis~\cite{Price}. The most
important result for our current purposes is the following theorem.
\begin{theorem}[\cite{str:plp}]
\label{layered} Let $q$ be a layered pattern. Then the
set of $q$-optimal $n$-permutations contains at least one layered
permutation.
\end{theorem}
Layered $\permpat$-optimal permutations have a simple recursive structure.
This comes from the fact, which we will use many times, that to form a
$\permpat$ pattern in a layered permutation one must take a single element
from some layer and $l$ elements from a subsequent layer
\begin{proposition} \label{fund}
Let $p$ be a layered $\permpat$-optimal $n$-permutation whose last layer is
of length $m$. Then the leftmost $k=n-m$ elements of $p$ form a
$\permpat$-optimal $k$-permutation.
\end{proposition}
\begin{proof} Let $D_k$ be the number of $\permpat$-copies of $p$ that are
disjoint from the last layer. The number of $\permpat$-copies of $p$ is
clearly $k {m \choose l}+D_k.$ So once $k$ is chosen, $p$ will have the
maximum number of copies only if $D_k$ is maximal.
\end{proof}

We point out that the proof of this proposition uses the fact that
$\permpat$ 
has only two layers, the first of which is a singleton.
Let $M_n=M_{n,\permpat}$.  Then the previous proposition implies that
\beq\label{Mn}
M_n=\max_{1\le k<n}\left( M_k+k {m \choose l}\right).
\eeq
The integer $k$ for which the right hand side attains its maximum will play
a 
crucial role throughout this paper. Therefore, we introduce specific
notation 
for it.
\begin{definition} For any positive integer $n$, let
$k_n=k_{n\permpat}$ be the
positive integer for which  $M_n=\max_k( M_k+ k {m \choose l})$ is
maximal. If there are several integers with this property,
then let $k_n$ be the largest among them.
\end{definition}
In other words, $k_n$ is the largest possible length of the
remaining permutation after removing the last layer of a $\permpat$-optimal
$n$-permutation $p$. When there is no danger of confusion, we will only
write 
$k$ to simplify notation. We will also always use $m=n-k$ to denote the
length of the last layer of $p$.
\section{ Construction of Permutations with a given number of copies
of $q=132$} 
We will first show that if $q=132$ then there are infinitely many
integers $n$ such that
$F_n$ does not have internal zeros. We will call such an integer, or its
corresponding sequence, {\it NIZ (no internal zero)}, and otherwise {\it
IZ}.  Our strategy is recursive: We will show that if $k_n$ is NIZ, then so
is $n$. As $k_n<n$, this will lead to an infinite  sequence of NIZ integers.
There is a problem, however. In order for this strategy to work, we must
ensure that given $k$, then there is an $n$ such that $k=k_n$. This is the
purpose of the following theorem which is in fact true for the general
pattern $q=\permpat$.

\begin{theorem} \label{cont}
For $k_n=k_{n,\permpat}$, the sequence $(k_n)_{n\ge1}$ diverges to
infinity and satisfies
$$
k_n\le k_{n+1}\le k_n+1
$$
for all $n\ge l+1$.  So, since $k_{l+1}=1$, for all positive
integers $k$ there is a positive integer $n$ so that $k_n=k$.
\end{theorem}

The next section is devoted to a proof of this theorem.
We suggest that the reader assume the result now and continue with this
section to preserve continuity.
We now consider the case $q=132$ which behaves differently from
$q=\permpat$ for $l\ge3$.  This is essentially due to the difference
between the patterns $q=12$ and $q=12\ldots l$ for $l\ge3$ as seen in
Proposition~\ref{12} and Theorem~\ref{12...l}.  First we note the
useful fact that
\beq                            \label{Mk}
M_{k,132}\ge{k-1\choose 2}
\eeq
which follows by considering the permutation
$1k(k-1)(k-2)\cdots 32$.
\begin{theorem} \label{nonzero}
For $q=132$ There are infinitely many NIZ integers.
\end{theorem}
\begin{proof}
It is easy to verify that $n=4$ is NIZ.
So, by Theorem~\ref{cont}, it suffices to show
that if $k_n\ge4$ is NIZ then
so is $n$.  To simplify notation in the two proofs which follow, we will
write $k$ for $k_{n,132}$, $M_n$ for $M_{n,132}$, and so forth.

Now given $c$ with $0\le c\le M_n = M_k+k{m\choose 2}$ we will
construct a permutation $p\in S_n$ having $c$ copies of 132.
Because of~\ree{Mk} and $k\ge4$  we have $M_k\geq k-1$.   So it is
possible to write $c$ (not necessarily uniquely) as $c=ks+t$ with
$0\le s\le {m\choose2}$ and $0\le t\le M_k$.  Since $k$ is NIZ, there
is a permutation $p'\in S_k$ with $c_{132}(p')=t$.  Also, by
Lemma~\ref{12/132}, there is a permutation in $S_m$ with no copies of 132
and
$s$ copies of 21.  Let $p''$ be the result of adding $k$ to every
element of that permutation.
Then, by construction, $p=p'p''\in S_n$ and $c_{132}(p)=ks+t=c$
as desired.
\end{proof}

One can modify the proof of the previous theorem to locate precisely
where the internal zeros could be for an IZ sequence.
We will need the fact (established by computer) that for $n\le 12$ the
only IZ integers were 6, 8, and 9, and that they all satisfied the
following result.
\begin{theorem} For any positive integer $n$, the sequence $F_{n,132}$ does
not
have internal zeros, except possibly for $c=M_{n,132}-1$ or
$c=M_{n,132}-2$, but not both.
\end{theorem}
\begin{proof}
We prove this theorem by induction on $n$. As previously remarked, it
is true if $n\leq 12$. Now suppose we know the statement for all integers
smaller than $n$, and prove it for $n$. If $n$ is NIZ, then
we are done.

If $n$ is IZ then, by the proof of Theorem~\ref{nonzero},
$k=k_n$ is IZ.  So $k\geq 6$ and we have
$M_k\ge k+2$ by~\ree{Mk}.  Now take $c$ with $0\le c\le M_n-3$
so that we can write $c=ks+t$ with
$0\le s\le {m\choose2}$ and $0\le t\le M_k-3$.  Since the portion of
$F_k$ up to $S_k(M_k-3)$ has no internal zeros by induction, we
can use the same technique as in the previous theorem to construct a
permutation $p$ with $c_{132}(p)=c$ for $c$ in the given range.
Furthermore, this construction shows that if $S_k(M_k-i)\neq 0$
for $i=1$ or $2$ then $S_n(M_n-i)\neq 0$.  This completes the proof.
\end{proof}

\section{The sequence $(k_n)_{n\ge l+1}$ for  $q=\permpat$}
For the rest of this paper, all invariants will refer to the
pattern $q=\permpat$ unless explicitly stated otherwise.

In order to prove Theorem~\ref{cont}, we first need a lemma about
the lengths of various parts of a $1(l+1)l\ldots 2$-optimal permutation $p$.
In all that follows,  we use the notation
\bea
b&=&\mbox{the length of the penultimate layer of $p$}\\
a&=&\mbox{the length of the permutation gotten by removing the last
        two layers of $p$}\\
 &=&n-m-b\\
 &=&k-b.
\eea
Also observe that the sequence $(M_n)_{n\ge l+1}$ is strictly increasing.
This 
is because when $n\ge l+1$, any layered $\permpat$-optimal
permutation $p\in S_n$  contains at least one copy of $\permpat$.  So
inserting $n+1$ in front of any layer contributing to the
$(l+1)l\ldots2$ portion of some copy results in a permutation with
more $\permpat$-patterns than $p$.  It follows from~\ree{Mn} that
$m\ge l$ for $n\ge l+1$, a fact that  will be useful in proving the
following result. 
\ble\label{le}
Let $q=1(l+1)l\ldots 2$, $k=k_{n,q}$, and $n\ge l+1$.
Then we have the following inequalities
\ben
\item[(i)] $b\le m$,
\item[(ii)] $a\le (m-l+1)/l$,
\item[(iii)] $k<n/l$, so in particular $k<m$,
\item[(iv)] $m\le l(n+1)/(l+1)$.
\een
\ele
\begin{proof} The basic idea behind all four of the inequalities is as
follows.  Let $p'$ be the permutation obtained from our $1(l+1)l\ldots 2$-
optimal
permutation $p$ by replacing its last two layers with a last layer of
length $m'$ and a next-to-last layer of length $b'$.  Then in passing
from $p$ to $p'$ we lose some $1(l+1)l\ldots 2$-patterns and gain some.
Since 
$p$ was optimal, the number lost must be at least as large as the number
gained.  And this inequality can be manipulated to give the one desired.

For the details, the following chart gives the relevant information to
describe $p'$ for each of the four inequalities. In the second case,
the last two layers of $p$ are
combined into one, so the value of $b'$ is irrelevant.
$$
\barr{c|c|c}
m'      &b'     &\mbox{number of gained $1(l+1)l\ldots 2$-patterns $\le$
                                number of lost $1(l+1)l\ldots 2$-patterns}\\
\hline \hline
b         &m      &m{b\choose l}\le b{m\choose l}\rule{0pt}{20pt}
\\[10pt]
\hline
b+m     &\mbox{---}&a\left({m+b\choose l}-{m\choose l}-{b\choose
l}\right)\le 
b{m\choose l}\rule{0pt}{20pt}\\[10pt]
\hline
m+1     &b-1    &(a+b-1){m\choose l-1}\le a{b-1\choose l-1}+{m\choose l}
\rule{0pt}{20pt}\\[10pt]
\hline
m-1     &b+1    &a{b\choose l-1}+{m-1\choose l}\le (a+b){m-1\choose l-1}
\rule{0pt}{20pt}
\earr
$$

Now (i) follows easily by cancelling $bm/l!$ from the inequality in
the first row of the table.

{From} the second line of the table, we have
$$
ab{m\choose l-1}\le
a\sum_{i=1}^{l-1} {b\choose i}{m\choose l-i}=
a\left({m+b\choose l}-{m\choose l}-{b\choose l}\right)\le b{m\choose l},
$$
and cancelling $b{m\choose l-1}$, which is not zero becase
$m\ge l$, gives us (ii).

To prove (iii) we induct on $n$.  If $n=l+1$, then we must have
$p=\permpat$, 
so $k=1<(l+2)/l=(n+1)/l$.  Now we assume $n>l+1$.

If $k<l+1$, then the leftmost $k$ elements of $p$ contain no copies of
$\permpat$, so we may replace them with any $k$-permutation and still have
$p$ optimal.  Therefore we may pick $b=1$ and $a=k-1$, and thus the second
row of the table shows
$$
\frac{l(k-1)}{m-l+1}{m\choose l}=(k-1){m\choose l-1}=(k-1)\left({m+1\choose
l}-{m\choose l}-{1\choose l}\right)\le {m\choose l},
$$
so $k\le (m+1)/l\le (n+1)/l$, as desired.

If $k\ge l+1$, recall that from Proposition~\ref{fund}, the leftmost $k=a+b$
elements of $p$ form a $\permpat$-optimal permutation, so we may, without
loss, choose $a$ maximal and thus assume that $a=k_{k}$.

{From} the third line of the chart, we have
$$
\frac{l(k-1)}{m-l+1}{m\choose l}=(a+b-1){m\choose l-1}\le
a{b-1\choose l-1}+{m\choose l}.
$$
Using (i) we get that
${b-1\choose l-1}\le {m-1\choose l-1}=\frac{l}{m}{m\choose l}$.
Substituting this in the previous
equation, cancelling ${m\choose l}$, and solving for $k$ gives
$$
k \le \frac{m+1}{l}+\frac{a(m-l+1)}{m}
$$
Since $k\ge l+1$, we have by induction that $a=k_k<k/l$.  Substituting
and solving for $k$ again and then cancelling $m+1$, we get
$k<\frac{m}{l-1}$.  A final substitution of $m=n-k$ results in (iii).

For (iv), notice that the last row of the table gives
\beq\label{m<eq}
{m-1\choose l}\le a{b\choose l-1}+{m-1\choose l}\le(a+b){m-1\choose
l-1}=(n-m)
{m-1\choose l-1}.
\eeq
so cancelling ${m-1\choose l-1}$ gives $n-m\ge(m-l)/l$, which
can be converted to the desired inequality.
\end{proof}

We now turn to the proof of Theorem~\ref{cont}.  First note that, by
Lemma~\ref{le} (iv), we have
\beq\label{k_lower_bound}
k=n-m\ge\frac{n-l}{l+1}.
\eeq
So $(k_n)_{n\ge1}$ clearly diverges to infinity.  For our next step, we
prove 
that $(k_n)_{n\ge1}$ is monotonically weakly increasing.  Let  $p_{n,i}$
denote an $n$-permutation whose last layer is of length $n-i$, and whose
leftmost $i$ elements form a $\permpat$-optimal $i$-permutation, and let $c_
{n,i}=c_{\permpat}(p_{n,i})$.  Clearly
$$
c_{n,i}=M_i+ i{n-i \choose l}.
$$

\begin{proposition} \label{kmon}
For $q=\permpat$ and all integers $n\ge l+1$, we have $k_n \le k_{n+1}$.
\end{proposition}
\begin{proof}
It suffices to show that $c_{n+1,k}>c_{n+1,i}$ for all $i< k$. This is
equivalent to showing that
\beq \label{kmon-1}
M_k + k{n-k+1\choose l} >M_i + i{n-i+1\choose l} .
\eeq
However, by definition of $k$, we know that for all $i<k$,
\beq \label{kmon-2}
M_k + k{n-k\choose l} \geq M_i + i{n-i\choose l} .
\eeq
Subtracting (\ref{kmon-2}) from (\ref{kmon-1}),  we are reduced to proving
$k{n-k\choose l-1} > i{n-i\choose l-1}$.
We will induct on $k-i$.  If $k-i=1$, then we would like to show that
$$
\frac{k(n-k-l+2)}{n-k+1}{n-k+1\choose l-1}=k{n-k\choose l-1}>(k-1){n-k+1
\choose l-1},
$$
so it suffices to show that $k<(n+1)/l$, which follows from Lemma~\ref{le}
(iii).

For $k-i>1$ we have, by induction, that $k{n-k\choose l-1} > (i+1){n-i-1
\choose l-1}$, so it suffices to show that
$$
\frac{(i+1)(n-i-l+1)}{(n-i)}{n-i\choose l-1}=(i+1){n-i-1\choose l-1}>i{n-
i\choose l-1},
$$
which simplifies to $(i+1)<(n+1)/l$, and this is is true because $i+1<k$.
\end{proof}

The proof of the upper bound on $k_{n+1}$ is a bit more involved but
follows the same general lines as the previous demonstration.  Note
that this will finish the proof of Theorem~\ref{cont}.

\begin{lemma} \label{0-1}
For $q=\permpat$ and all integers $n\ge l+1$, we have $k_n\le k_{n+1}\le
k_n+1$.
\end{lemma}
\begin{proof}
Induct on $n$.  The lemma is true for $n=l+1$ since $k_{l+1}=k_{l+2}=1$.
Suppose  the lemma is true for integers
smaller than or equal to $n$, and prove it for $n+1$.  For simplicity, let
$k=k_n$, $m=n-k$, and $c_i=c_{n+1,i}$. Since we have already proved the
lower 
bound, it suffices to show that
\beq\label{ci}
c_i\ge c_{i+1}\mbox{ for } k+1\le i<\flf{n+1}{l}
\mbox{ with strict inequality for } i=k+1.
\eeq
Note that we do not have to consider $i\ge\fl{(n+1)/l}$
because of Lemma~\ref{le} (iii).

We prove~\ree{ci} by induction on $i$.  For the base case, $i=k+1$,
we wish to show
\beq                    \label{Mk+1/2}
M_{k+1}+(k+1){m\choose l}> M_{k+2}+(k+2){m-1\choose l}.
\eeq
But since $p_{n,k}$ is optimal by assumption, we have
\beq                    \label{Mk/k+1}
M_k+k{m\choose l}> M_{k+1}+(k+1){m-1\choose l}.
\eeq
Subtracting~\ree{Mk/k+1} from~\ree{Mk+1/2} and rearranging terms, it
suffices to prove
\beq                            \label{m-1}
{m-1\choose l-1}\ge (M_{k+2}-M_{k+1})-(M_{k+1}-M_k).
\eeq

First, if $k<l+1$, then~\ree{m-1} is easy to verify using
Lemma~\ref{le} (iii) and the values $M_{l+2}=l+1$, $M_{l+1}=1$, and
$M_k=0$ for $k\le l$.
Therefore we may assume that $k\ge l+1$.  Let $p'\in S_k$, $p''\in S_{k+1}$,
and $p'''\in S_{k+2}$ be layered $\permpat$-optimal permutations having last
layer lengths $m'$, $m''$, and $m'''$, respectively, as short as possible.
Also let $k'=k-m'$, $k''=k+1-m''$, and $k'''=k+2-m'''$.  We would like to be
able to assume the lemma holds for these permutations, and thus we would
like 
to have $k+2\le n$.  But by Lemma~\ref{le} (iii) we have $k+2<n/2+2\le
n$ if $n\ge4$.  Since $n\ge l+1$ this holds for $l\ge3$ and the case
$l=2, n=3$ is easy to check directly.
Therefore we may assume that $p'$, $p''$, and $p'''$ all satisfy the lemma.

If $m''=m'+1$ then let $x$ be the largest element in the last layer of $p''$
(namely $x=k+1$).  Otherwise, $m''=m'$ and removing the last layer of both
$p'$ and $p''$ leaves permutations in $S_{k-m'}$ and $S_{k-m'+1}$,
respectively.  So we can iterate this process until we find the single layer
where $p'$ and $p''$ have different lengths (those lengths must differ by 1)
and let $x$ be the largest element in that layer of $p''$.   Similarly we
can 
find the element $y$ which is largest in the unique layer were $p''$ and
$p'''$ have different lengths.

Now let
\bea
r&=&
\mbox{the number of $\permpat$-patterns in $p'''$ containing neither $x$ nor
$y$,}\\
s&=&
\mbox{the number of $\permpat$-patterns in $p'''$ containing $x$ but not
$y$,}
\\
t&=&
\mbox{the number of $\permpat$-patterns in $p'''$ containing $y$ but not
$x$, 
and}\\
u&=&\mbox{the number of $\permpat$-patterns in $p'''$ containing both $x$
and 
$y$.}
\eea
Note that there is a bijection between the $\permpat$-patterns of $p'''$ not
containing $y$ and the $\permpat$-patterns of $p''$.  A similar statement
holds for $p''$ and $p'$.  So
$$
M_k=r,\quad M_{k+1}=r+s,\quad M_{k+2}=r+s+t+u.
$$
Note also that $s\ge t$ because increasing the length of the layer of
$x$ results in the most number of $\permpat$-patterns being added to $p'$.
It follows that  $(M_{k+2}-M_{k+1})-(M_{k+1}-M_k)=t+u-s\le u$.

By Lemma~\ref{le} (iii), $k<m$, so to obtain~\ree{m-1}
it suffices to show that $u\le {k\choose l-1}$.  But ${k\choose l-1}$
is the total number of subsequences of $p'''$ having length $l+1$ and
containing $x$ and $y$.  So the inequality follows.

The proof of the induction step is similar.  Assume that~\ree{ci} is true
for 
$i-1$ so that
\beq\label{ci/i+1}
M_{i-1}+(i-1){r+1\choose l}\ge M_i+i{r\choose l}.
\eeq
where $r=n+1-i$.  We wish to prove
\beq
M_i+i{r\choose l}\ge M_{i+1}+(i+1){r-1\choose l}.
\eeq
Subtracting as usual and simplifying, we need to show
$$
2{r-1\choose l-1}-(i-1){r-1\choose l-2}\ge (M_{i+1}-M_i)-(M_i-M_{i-1}).
$$
Proceeding exactly as in the base case, we will be done if we can show
that
$$
\frac{2r-l-il+i+1}{r-l+1}{r-1\choose l-1}=2{r-1\choose l-1}-(i-1){r-1\choose
l-2}\ge{i-1\choose l-1}.
$$
Because $i<\flf{n+1}{l}$ we have $r\ge i$, so it suffices to show that
$$
\frac{2r-l-il+i+1}{r-l+1}{r-1\choose l-1}\ge 1.
$$
This simplifies to showing that $i\le (r+i)/l=(n+1)/l$, and this is
guaranteed by our choice of $i$.
\end{proof}

The following lemma contains two inequalities essentially shown in the proof
of
Lemma~\ref{0-1} which we will need to use again.

\begin{lemma} \label{Mkconcave}
If $q=\permpat$ then
$0\le (M_{i+2}-M_{i+1})-(M_{i+1}-M_i)\le {i\choose l-1}$.
\end{lemma}
\begin{proof}
For the upper bound, recall that ${i\choose l-1}$ is the total number
of subsequences of $p'''$ of length $l+1$ containing $x$ and $y$ while
the double difference just counts those subsequences corresponding to
the pattern $q=\permpat$.
For the lower bound, we showed that
$$
(M_{i+2}-M_{i+1})-(M_{i+1}-M_i) = t + u - s.
$$
Recall that $t+u$ is the total contribution of $y$ in $p'''$, and $s$ is the
total contribution of $x$ in $p''$.  Therefore $t+u-s \ge 0$, as otherwise
one could create a permutation with more $\permpat$-patterns than $p'''$ by
inserting a new element in the same layer as $x$
\end{proof}
\section{The sequence $(c_{n,i})_{i=1}^{n-1}$ for $q=\permpat$}

Now that we have completed the proof of Theorem~\ref{cont}, we turn our
attention to the tools which will enable us to show that there are
infinitely 
many IZ integers.  As before, all invariants are for $q=\permpat$
unless otherwise stated.

For $l=2$, we will need the following lemma.

\begin{lemma} \label{1storderbounds}
For all $n$, we have $M_{n+1,132} - M_{n,132} \leq 5n^2/16$.
\end{lemma}
\begin{proof}
Let $k = k_n$.  We induct on $n$.  It is easy to check the base cases
$n=1,2$.  Note that by Theorem~\ref{cont}, either $k_{n+1} = k$ or
$k_{n+1} = k + 1$.  If $k_{n+1} = k$, then we have
$$
M_{n+1} - M_n = nk - k^2
$$
and maximizing this as a function of $k$ gives
$$
M_{n+1} - M_n \leq \frac{n^2}{4}\leq \frac{5n^2}{16}.
$$
If $k_{n+1} = k + 1$, then we have
$$
M_{n+1} - M_n = M_{k+1} - M_k + {n-k\choose 2}.
$$
By induction, we have $M_{k+1} - M_k \leq \frac{5k^2}{16}$, and thus we have
that
$$
M_{n+1} - M_n \leq \frac{13k^2}{16} + \frac{n^2 + k - n - 2kn}{2}.
$$
By Lemma~\ref{le} (iii) and (iv), this function is to be maximized on
the interval $[(n-2)/3,n/2]$ and for $n\ge3$ this maximum occurs at
$k=(n-2)/3$.  So
$$
M_{n+1} - M_n \leq \frac{37n^2-4n+4}{144} \leq \frac{5n^2}{16},
$$
as desired.
\end{proof}

\begin{definition} For $q=\permpat$ and any positive integer $n$, let
$l_n$ be the least 
integer greater than $k_n$ such that $c_{n,i} \le c_{n,i+1}$. If there is no
integer with this property, let $l_n=n-1$.
\end{definition}

Do not confuse $l_n$, which will always be subscripted, with the
length-related parameter $l$ , which will never be.
Our next result shows that the sequence $(c_{n,i})_{i=1}^{n-1}$ is
``bimodal''
with a maximum at $i=k_n$ and a minimum at $i=l_n$.


\begin{theorem} \label{bimodal}
For $q=\permpat$ and all positive integers $n$, we have the following three
results about
the shape of $(c_{n,i})_{i=1}^{n-1}$
\ben
\item[(i)] $c_{n,i} \le c_{n,i+1}$ for all $i < k_n$,
\item[(ii)] $c_{n,i} \ge c_{n,i+1}$ for all $k_n \le i < l_n$,
\item[(iii)] $c_{n,i} \le c_{n,i+1}$ for all $i\ge l_n$.
\een
\end{theorem}
\begin{proof}
For (i) we induct on $n$.  The claim is true trivially for $n<l+1$ since
then $c_{n,i}=0$ for all $i$, so we will assume $n\ge l+1$.  
If $i=k_n-1$ then the claim is true by definition.  If $i<k_n-1$ 
then $i<k_{n-1}$ by Theorem~\ref{cont} and we are able to apply induction.
We would like to show that
$$
M_i + i{n-i\choose l}\le  M_{i+1} + (i+1){n-i-1\choose l} 
$$
and we know by induction that
$$
M_i + i{n-i-1\choose l}\le  M_{i+1} + (i+1){n-i-2\choose l}. 
$$
Subtracting as usual, we are reduced to showing that 
$i{n-i-1\choose l-1}\le  (i+1){n-i-2\choose l-1}$.  This
further reduces to $i \le (n-l)/l$ which is true by
Lemma~\ref{le} (iii) and the fact that $i<k_n-1$.

Statement (ii) is implied by the definition of $l_n$, so we are left with
(iii).   By the definition of $l_n$ we have that
$c_{n,l_n} \le c_{n,l_n+1}$, so it suffices to show that 
for all $i\ge l_n$, if $c_{n,i}\le c_{n,i+1}$ then $c_{n,i+1}\le c_{n,i+2}$.  
Subtracting in the usual way, we are reduced to showing that
\beq \label{bimodal3DD}
(M_{i+2}-M_{i+1})-(M_{i+1}-M_i) \ge \frac{2n-2l-i(l+1)}{l-1}{n-i-2\choose l-
2}.
\eeq
Since we know that $(M_{i+2}-M_{i+1})-(M_{i+1}-M_i) \ge 0$ by Lemma~\ref
{Mkconcave}, our approach will be to show that $2n-2l-i(l+1) \le 0$ for $i \ge
l_n$ by showing that
\beq\label{ln_lower}
l_n \ge (2n-2l)/(l+1).
\eeq

Before we prove (\ref{ln_lower}), we will need the following two facts.
\begin{center}
$l_n\ge n/l$ and  $l_n\ge l_{n-1}$.
\end{center}
The first fact follows from our proof of Lemma~\ref{0-1}, in which we showed
that $c_{n,i}\ge c_{n,i+1}$ for $k_n\le i<\fl{n/l}$.
So to prove the second fact, it suffices to show
that  $c_{n-1,i}>c_{n-1,i+1}$ implies $c_{n,i}> c_{n,i+1}$
for $i\ge n/l$.  This is proved in exactly the same way as (i) with 
all the inequalities reversed.

Now we are ready to prove (\ref{ln_lower}).  First we tackle the case
where $l\ge3$ by induction. 
If $n\le3$ then $(2n-2l)/(l+1)\le0$ and we are done.  So suppose $n\ge4$.
If $l_{n-2}>(n-1)/2$, then since $l_n\ge l_{n-2}$ and $l\ge3$ we have
$l_n>(2n-2l)/(l+1)$ as
desired.  Hence we may assume that $(n-2)/l\le l_{n-2}\le (n-1)/2$.  In this
case we claim that $l_n\ge l_{n-2}+1$, which will imply (\ref{ln_lower})
by induction.

Let $i=l_{n-2}$.  We want to show that
$$
M_i+i{n-i\choose l} > M_{i+1}+(i+1){n-i-1\choose l},
$$
and we have
$$
M_{i-1}+(i-1){n-i-1\choose l} > M_i + i{n-i-2\choose l}.
$$
Subtracting, it suffices to show that
$$
(M_{i+1}-M_i)-(M_i-M_{i-1})\le i{n-i-2\choose l-2}.
$$
By Lemma~\ref{Mkconcave}, $(M_{i+1}-M_i)-(M_i-M_{i-1})\le {i-1\choose l-1}$,
so it suffices to show that
\beq\label{l>2ln-2}
{i-1\choose l-1}\le \frac{il-i}{n-i-l}{n-i-2\choose l-1}.
\eeq
Since $i\ge (n-2)/l$, we have that $(il-i)/(n-i-l)\ge
1$, and since $i\le (n-1)/2$, we have that $n-i-2\ge i-1$, so
(\ref{l>2ln-2}) is true, and thus (\ref{ln_lower}) holds.

For the case where $l=2$, we examine the quadratics
$$
d_i(n) = \frac{1}{2}n^2 - \left(2i+\frac{3}{2}\right)n + \left(M_{i+1}-
M_i+\frac{3}{2}i^2+\frac{5}{2}i+1\right),
$$
which agree with $c_{n,i+1}-c_{n,i}$, wherever both $c_{n,i+1}$ and
$c_{n,i}$ 
are defined.  We will also need to refer to the roots of $d_i(n)$, which
occur at
$$
r_i = 2i + \frac{3}{2} - \sqrt{i^2+i+\frac{1}{4}-2\left(M_{i+1}-
M_i\right)},\mbox{ and}
$$
$$
s_i = 2i + \frac{3}{2} + \sqrt{i^2+i+\frac{1}{4}-2\left(M_{i+1}-M_i\right)}.
$$
Lemma~\ref{1storderbounds} gives us that
\beq \label{three-eighths}
r_i < (2 - \sqrt{3/8})i + 3/2,
\eeq
so $r_i$ and $s_i$ are real numbers and for $i>13$,
$r_i<3i/2$.  These roots are important in our situation for the following
reasons:

\beq \label{di<0}
d_i(n) < 0 \mbox{ if and only if } r_i < n < s_i,
\eeq
\beq \label{n>si}
n \ge s_i \mbox{ if and only if } i \le k_n, \mbox{ and}
\eeq
\beq \label{n<=rln}
n\le r_{l_n}.
\eeq

Statement (\ref{di<0}) is easily verified.  Assume to the contrary that the
forward direction of (\ref{n>si}) is not true, and thus $n \ge s_i$ but $i >
k_n$.  Let $n'$ be such that $k_{n'} = i$.  By Proposition~\ref{kmon}, we
have that $n' > n \ge s_i$, and thus $d_i(n') \ge 0$ by (\ref{di<0}).
However because $i=k_{n'}$, we have that $d_i(n') < 0$, a contradiction.  To
prove the reverse direction of (\ref{n>si}), notice that if $i\le k_n$ then
by (i) and the definition of $k_n$, we must have that $d_i(n) \ge 0$.
Therefore by (\ref{di<0}), either $n\ge s_i$ (as we would like) or $n\le
r_i$, and by (\ref{three-eighths}), it cannot be the case that $n \le r_i$,
as that would imply that $n \le r_i < 3i/2 < 3k_n/2$ if $i>13$,
contradicting Lemma~\ref{le} (iii).  

To prove (\ref{n<=rln}), note that by (\ref{di<0}) we cannot have
$r_{l_n}<n<s_{l_n}$ as then we would have $d_{l_n}(n)<0$, contradicting the
definition of $l_n$.  Also, we cannot have $n\ge s_{l_n}$ as then we would
have $l_n\le k_n$ by (\ref{n>si}), again contradicting the definition of
$l_n$.  Hence we must have (\ref{n<=rln}).

With these tools, (\ref{ln_lower}) is easy to prove; we have $n\le
r_{l_n}<3l_n/2$ for $l_n>13$, and thus $l_n>2n/3$, as desired.  It is easily
checked that $l_n>2n/3$ for $l_n\le 13$.
\end{proof}

We will depend on the following lemma to find integers $n$ with an internal
zero at $M_n-1$.

\begin{lemma}\label{cnidiff}
For $q=\permpat$, $l\ge 2$ and all $n\ge 2l+2$, if $k_{n-2} = k-1$ and
$k_{n-1} = k$, then
$M_n - c_{n,i} > 1$ for all $i\ne k$, so in particular, $k_n=k$.
\end{lemma}
\begin{proof}
By Theorem~\ref{bimodal} it suffices to show the following inequalities:
\beq    \label{cnidiffk-1}
c_{n,k} - c_{n,k-1} > 1,
\eeq
\beq    \label{cnidiffk+1}
c_{n,k} - c_{n,k+1} > 1,
\eeq
and
\beq    \label{cnidiffn-1}
c_{n,k} - c_{n,n-1} > 1.
\eeq

Statement (\ref{cnidiffn-1}) is clear for $n \ge 2l+2$ because $c_{n,n-1} =
M_{n-1}$, $(M_i - M_{i-1}) \ge (M_{i-1}-M_{i-2})$ for all $i$ by
Lemma~\ref{Mkconcave}, and $M_{l+2}-M_{l+1} = l$.

We prove (\ref{cnidiffk+1}) by induction on $n$.  First, if $k<l$, then
$M_{k-1}=M_k=M_{k+1}=0$, so it suffices to show that
$$
k{n-k\choose l}>(k+1){n-k-1\choose l}+1,
$$
and since $k_{n-2}=k-1$, we have
$$
(k-1){n-k-1\choose l}>k{n-k-2\choose l}.
$$
Subtracting that latter from the former, it suffices to show that
$$
1\le k{n-k-2\choose l-2}.
$$
So we're done in this case since $n-k\ge l$ which follows from
$n\ge2l+2$ and $k<l$.

Now assume that $k\ge l$, so we may prove
(\ref{cnidiffk+1}) by showing the stronger statement that
$$
c_{n,k} - c_{n,k+1} > {k-2\choose l-2},
$$
and thus we would like to show that
$$
M_k + k{n-k\choose l} > M_{k+1}+(k+1){n-k-1\choose l} + {k-2\choose l-2},
$$
and as $k_{n-2} = k-1$, we have
$$
M_{k-1}+(k-1){n-k-1\choose l} > M_k + k{n-k-2\choose l}.
$$
Subtracting as usual, we are reduced to showing
$$
(M_{k+1}-M_k)-(M_k-M_{k-1})\le k {n-k-2\choose l-2}-{k-2\choose l-2}.
$$
By Lemma~\ref{le} (iii)
$$
k {n-k-2\choose l-2}-{k-2\choose l-2}>
k {k-2\choose l-2}-{k-2\choose l-2}=(l-1){k-1\choose l-1}.
$$
The upper bound in Lemma~\ref{Mkconcave} now completes the proof of
(\ref{cnidiffk+1}).

To prove (\ref{cnidiffk-1}), we want to show
$$
M_k + k{n-k\choose l} > M_{k-1}+(k-1){n-k+1\choose l} + 1,
$$
and we are given
$$
M_{k}+k{n-k-1\choose l} \ge M_{k-1} + (k-1){n-k\choose l}.
$$
Subtracting as usual, we are reduced to showing that
$$
k{n-k-1\choose l-1}>(k-1){n-k\choose l-1} + 1.
$$
Cancelling ${n-k-1\choose l-1}$ and simplifying, it suffices to show that
\beq  \label{branch}
n>lk + \frac{n-k-l+1}{{n-k-1\choose l-1}}.
\eeq
By Lemma~\ref{le} (iii), $n\ge lk+1$, so it suffices to show
that
$$
n-k-l+1<{n-k-1\choose l-1},
$$
which is true for $l\ge3$.  For $l=2$, note that proving~\ree{branch}
reduces to showing $n>2k+1$ which we will prove by induction on $n$.
Checking the base cases $n=6,7$ is easy.  Also note that~\ree{branch}
holds for $l=2$ if we make the strict inequality
weak, so we still can conclude the $k_n=k$ part of the Lemma.
There are now two cases.  If $k_{n-1}=k_n=k$ then by induction
$n>n-1>2k+1$.  By Theorem~\ref{cont} and the part of the Lemma that
we've already proved, the only other possibility is $k_{n-1}=k-1$ and
$k_{n-2}=k-1$.  But then $n-2>2(k-1)+1$ which is equivalent to the
desired inequality.
\end{proof}

\section{The poset connection}

There is an intimate connection between partially ordered sets, called
posets 
for short, and permutations.  Using this connection, we will provide
characterizations of all $n$-permutations $p$ which have $c_{1(l+1)l\dots
2}(p)
\ge M_{n,1(l+1)l\dots 2}-1$ for $l\ge 2$.  This will provide us with the
tools 
we need to show that there are an infinite number of IZ sequences for each
of 
these patterns. Any necessary definitions from the theory of posets that are
not given here will be found in Stanley's text~\cite{sta:ec1}.

If $P$ is a poset such that any two distinct elements of $P$ are
incomparable 
we say that $P$ is an {\it antichain\/}.  Since there is a unique unlabelled
antichain on $n$ elements, we denote this poset  by $\achain_n$.

Given posets $P$ and $Q$, the {\it ordinal sum\/} of $P$ and $Q$, denoted
$P\oplus Q$, is the unique poset on the elements $P\cup Q$ where $x\le y$ in
$P\oplus Q$ if either
\ben
\item[(i)] $x,y\in P$ with $x\le y$,
\item[(ii)] $x,y\in Q$ with $x\le y$, or
\item[(iii)] $x\in P$ and $y\in Q$.
\een
A poset $P$ is {\it layered\/} if
it is an ordinal sum of antichains, i.e. if $P=\achain_{p_1}\oplus\achain_
{p_2}\oplus\cdots\oplus\achain_{p_k}$ for some $p_1,\dots,p_k$.  To
introduce 
a related notion, let $\max P$ denote the set of maximal elements of $P$ and
$\minusmax{P}=P\setminus(\max P)$.  Then $P$ is {\it LOT (layered on top)}
if 
$P=\minusmax{P}\oplus\max P$.  Note that if $P$ is layered then it is LOT,
but not conversely.

If $p=p_1p_2\ldots p_n$ is a permutation, then the corresponding poset $P_p$
has elements $p_1,p_2,\ldots,p_n$ with partial order $p_i<p_j$ if
$(p_i,p_j)$ 
is a noninversion in $p$.  So, for example, $P_{12\ldots n}$ is a chain, $P_
{n\ldots 21}=\achain_n$ and $P_{1(l+1)l\dots 2}=\pat$.  Clearly not every
poset is of the form $P_p$ for some $p$.  In fact, the $P_p$ are exactly the
posets of dimension at most 2, being the intersection of the total orders
$1<2<\cdots<n$ and $p_1<p_2<\cdots<p_n$.

Given posets $P$ and $Q$ let
$$
\mbox{$c_Q(P) =$ the number of induced subposets of $P$ isomorphic to $Q$}.
$$
Now given permutations $p,q$ with corresponding posets $P=P_p, Q=P_q$,
we have $c_q(p)\le c_Q(P)$ since the elements of each copy of
$q$ in $p$  form a subposet of $P$ isomorphic to $Q$.

If $S \subseteq P$ then let
\bea
c_Q(P;S) &=&\mbox{the number of induced $Q'\subseteq P$ with
    $Q'\iso Q$ and $S \cap Q' \neq \emptyset$},\\
c_Q(P;\newnot S) &=&\mbox{the number of induced $Q'\subseteq P$ with
    $Q'\iso Q$ and $S \cap Q' = \emptyset$}.
\eea
We will freely combine these notations and eliminate the subscript when
talking about a fixed poset $Q$.  We will also abbreviate $c_Q(P;\{x\})$ to
$c_Q(P;x)$ and $c_Q(P;\newnot \{x\})$ to $c_Q(P; \newnot x)$.

As with permutations, for any non-negative integer $n$ we will let $M_{n,Q}
=\max\{c_Q(P) : |P| = n\}$.  We will say a poset $P$ is $Q$-optimal if
$c_Q(P)
=M_{|P|,Q}$.

Stromquist proved Theorem~\ref{layered} by first demonstrating the following
stronger result.

\begin{theorem}[\cite{str:plp}]
\label{strom}
If $Q$ is a LOT pattern, then there is some $Q$-optimal LOT poset $P$.  The
same holds  with ``LOT'' replaced by ``layered.''
\end{theorem}

To show that the sequences of the patterns $1(l+1)l\dots 2$, for
$l\ge 2$,
have infinitely many IZ integers, we will need to know more about
$\pat$-optimal 
posets. The best possible case would be if all (sufficiently large)
$\pat$-optimal 
posets were layered.  This is true for the pattern
$P_{132}=\achain_1\oplus\achain_2$, but not in general.  For example, it can
be computed that $P_{231}\oplus\achain_{8}$ is
$\achain_1\oplus\achain_3$-optimal, but $P_{231}\oplus\achain_{8}$ is not
layered.  Fortunately, we are able to show that all $\pat$-optimal posets
are of the following slightly more general form.

\begin{definition}
We say $P=P_1\oplus P_2$ is an $l$-decomposition of $P$ if $P_2$ is layered
and for all $A\subseteq P$ with $A\cong \pat$ we have
$|A\cap P_1|\le 1$.
\end{definition}

The first part of this section concerns the proof of the following theorem.

\begin{theorem}\label{strongstrom}
If $P$ is an $\pat$-optimal poset then $P$ has an $l$-decomposition.
\end{theorem}

After this proof we will investigate `almost' $\pat$-optimal posets, that
is, 
posets $P$ with $c_{\pat}(P)=M_{|P|,\pat}-1$.

If $q$ and $p$ are permutations, it is generally not the case that
$c_{P_q}(P_p)=c_q(p)$. For example, $P_{231}\cong P_{312}$ and thus
$c_{P_{231}}(P_{312})=1$, but $c_{231}(312)=0$.  However, there is an
important case in which we do get equality.

\begin{lemma}\label{permposetequiv}
If $q$ and $p$ are permutations then $c_q(p)\le c_{P_q}(P_p)$.  Furthermore,
if either $q$ or $p$ is layered then $c_q(p)=c_{P_q}(P_p)$.
\end{lemma}
\begin{proof}
The inequality follows from the fact that each copy of $q$ in $p$
gives rise to a copy of $P_q$ in $P_p$.  For the equality, if $q$ is
layered then it is the unique permutation giving rise to the poset
$P_q$.  So every copy of $P_q$ in $P_p$ corresponds to a copy of $q$
in $p$ and we are done.  The only other case we need to consider is if
$p$ is layered and $q$ is not.  But then both sides of the equality
are zero.
\end{proof}

This lemma and the preceeding theorems imply several important features
about 
the connection between pattern matching in posets and permutations.
Given any pattern $q$, the first statement in
Lemma~\ref{permposetequiv} implies that $M_{n,q}\le M_{n,P_q}$ for all
$n$.
If $q$ is layered, then by Theorem~\ref{strom} there is a layered
$P_q$-optimal poset $P=\achain_{p_1}\oplus\achain_{p_2}
\oplus\cdots\oplus\achain_{p_k}$ for some positive integers $p_1,\dots,p_k$.
It follows that there is a layered permutation $p$ such that $P_p\iso P$,
namely $p$ is the permutation whose layer lengths from left to right are
$p_1,\dots,p_k$.  By the preceeding lemma, $c_q(p)=c_{P_q}(P)$, so
$M_{|P|,q}=M_{|P|,P_q}$.

\begin{lemma}\label{Mmonotone}
For all patterns $Q$, the sequence $(M_{n,Q})_{n\ge|Q|}$ is positive
and strictly increasing.
\end{lemma}
\begin{proof}
We will write $M_n$ for $M_{n,Q}$ and $c(P)$ for $c(P)$.  Given
$n\ge|Q|$, it is easy to construct a poset $P$ with $c(P)>0$.  So let
$P$ be a $Q$-optimal poset.
Now there must be some $x\in P$ with $c(P;x)>0$.  Now adjoin an
element  $y$  to $P$ to form a poset
$P'$ with $a< b$ in $P'$ if
either 
\ben
\item[(i)] $a,b\in P$ with $a< b$,
\item[(ii)] $a=y$, $b\in P$ with $x< b$, or
\item[(iii)] $b=y$, $a\in P$ with $a< x$.
\een
Then
$$
c(P')=c(P';\newnot y)+c(P';y)=c(P)+c(P';y)\ge c(P)+c(P;x)>c(P)
$$
so $M_{n+1}\ge c(P')>c(P)=M_{n}$.
\end{proof}

We now begin the proof of Theorem~\ref{strongstrom} by making a few
definitions.  If $P$ is a poset and $x\in P$ then the {\it open down-set\/}
generated by $x$ is
$$
P_{<x}=\{y\in P\ :\ y<x\}.
$$
If $x,y\in\max P$ then let $P^{x\ra y}$ be the unique poset on the same set
of elements which satisfies
$$
\mbox{$P_{<z}^{x\ra y} = P_{<z}$ for $z\neq x$ and $P_{<x}^{x\ra y} = P_{<y}
$}.
$$
Note that $P-x=P^{x \rightarrow y}-x$.
The following lemma is essentially in Stromquist~\cite{str:plp}, but is not
explicitly proved there.  So we will provide a demonstration.
\begin{lemma}
\label{stromLOT}
Let $Q$ be a LOT pattern and $P$ be any poset with $x,y \in \max P$.  Then
$$
c_Q(P^{x \rightarrow y}) \geq c_Q(P) + c_Q(P;y) - c_Q(P;x).
$$
\end{lemma}
\begin{proof}
As before, we write $c(P)$ for $c_Q(P)$.
Since
\begin{eqnarray*}
c(P)              &=& c(P;\newnot x) + c(P;y,\newnot x) + c
(P;x,y),\mbox{ and}\\
c(P^{x \rightarrow y})    &=& c(P^{x \rightarrow y};\newnot x) +
c(P^{x\rightarrow y};x,\newnot y) + c(P^{x \rightarrow y};x,y),
\end{eqnarray*}
it is enough to show that
\begin{eqnarray}
c(P^{x \rightarrow y}; \newnot y)
&=&
c(P;\newnot y),
\label{stromeqn1}   \\
c(P^{x \rightarrow y}; x, \newnot y)
&\geq&
c(P;x,\newnot y) + c(P;y) - c(P;x),
\label{stromeqn2} \\
c(P^{x \rightarrow y}; x,y)
&\geq&
c(P;x,y).
\label{stromeqn3}
\end{eqnarray}

First, (\ref{stromeqn1}) is clear since $P$ and $P^{x \rightarrow y}$ agree
on all subsets not including $x$.

Next, notice that
$$
c(P;x,\newnot y) + c(P;y) - c(P;x) = c(P;y,\newnot x),
$$
and thus to prove (\ref{stromeqn2}), it suffices to show that $c(P^{x
\rightarrow y};x,\newnot y) \geq c(P;y,\newnot x)$, but this is easy.
Let $A \subseteq P$ with $y \in A$, $x \notin A$, and $(A, \leq) \cong Q$.
Then $A' = A \cup \{ x \} - y$ is an occurance of $Q$ in $P^{x \rightarrow
y}$, i.e., $(A', \leq_{P^{x \rightarrow y}}) \cong Q$, so (\ref{stromeqn2})
is
proved.

Finally, to prove (\ref{stromeqn3}), let $A \subseteq P$ be an occurance of
$Q$ in $P$ which contains $x$ and $y$, i.e., $(A,\leq_P)\cong Q$.  Then
we have that $(A, \leq_{P^{x \rightarrow y}}) \cong Q$ as well.  This is
because $A_{<x}=A_{<y}$ in $P$ since $x,y$ are maximal and $Q$ is LOT.
So $A$ forms an occurance of $Q$ in $P^{x \rightarrow y}$, and thus
(\ref{stromeqn3}) is proven.
\end{proof}

For the rest of this section, let $Q_l=\pat$,
$c(P)=c_{Q_l}(P)$ and
$M_n=M_{n,Q_l}$.

\begin{lemma}\label{minusa}
Let $P$ be a poset such that
$|P|>l\ge2$. If for some $a\ge 0$ and $x\in\max P$ we
have
\ben
\item[(a)] $P-x$ is LOT,
\item[(b)] $\packpat(P)=M_{|P|,Q_l}-a$, and
\item[(c)] $\packpat(P;x)=\packpat(P;y)-a$ for all
$y\in\max P\setminus x$,
\een
then $P$ is LOT (and thus $a$ is actually $0$).
\end{lemma}
\begin{proof}
Choose $y\in\max P$ with $y\ne x$, let $m=|\max P|$ and
$k=|\minusmax{P}|=|P|-m$.  First consider what happens when $m<l$.
Then (a) implies that $C(P;y)=0$ for all $y\in\max P\setminus x$.
This forces $c(P;x)=a=0$ by (c).  Now (b) yields
$c(\bar{P})=c(P)=M_{|P|}$, contradicting Lemma~\ref{Mmonotone}.  So we
may assume $m\ge l$.

Note that $\packpatabbrev(P)=\packpatabbrev(P-x)+\packpatabbrev(P;x)
$, and since 
$\packpatabbrev(P;x)=\packpatabbrev(P;y)-a=\packpatabbrev(P;x,y)
+\packpatabbrev(P;\newnot
x,y)-a=\packpatabbrev(P;x,y)+\packpatabbrev(P-x;y)-
a$ we get that
$$
\packpatabbrev(P) = \packpatabbrev(P-x)+\packpatabbrev(P-x;y)+\packpatabbrev
(P;x,y)-a.
$$
Furthermore, since $P-x$ is LOT we get that
$$
\packpatabbrev(P-x) = \packpatabbrev(\minusmax{P}) + {m-1\choose l}k,
$$
and
$$
\packpatabbrev(P-x;y) = {m-2\choose l-1}k.
$$
Also, since $P_{<x}\subseteq\minusmax{P}=P_{<y}$, we have that
$$
\packpatabbrev(P;x,y) = {m-2\choose l-2}|P_{<x}|,\mbox{ so}
$$
\beq\label{fP}
\packpatabbrev(P) = \packpatabbrev(\minusmax{P}) + \left( {m-1\choose
l}k+{m-2
\choose l-1}k + {m-2\choose l-2}|P_{<x}|\right) - a.
\eeq

Furthermore, since $P-x$ is LOT, $\pxy$ is LOT, so we have
$$
\packpatabbrev(\pxy) = \packpatabbrev(\minusmax{P}) + {m\choose l}k.
$$
Therefore
\begin{eqnarray}
\packpatabbrev(\pxy) - \packpatabbrev(P) &=& a + \left( \left( {m\choose l}
- 
{m-1\choose l} - {m-2\choose l-1} \right) k - {m-2\choose
l-2}|P_{<x}|\right) \nonumber
\\
&=& a + {m-2\choose l-2}\left(k - |P_{<x}|\right). \label{for7.9}
\end{eqnarray}
Furthermore, by Lemma~\ref{stromLOT}  and assumptions (b) and (c) we
have that $\packpatabbrev(\pxy)\ge M_{|P|}$.  So we must have
$\packpatabbrev(\pxy)=M_{|P|}$ and, by (b) again,
$\packpatabbrev(\pxy)-\packpatabbrev(P)=a$.  It follows that
${m-2\choose l-2}\left(k-|P_{<x}|\right) =0$.  Therefore since $m\ge
l$ we have ${m-2\choose l-2}>0$ and so $k=|P_{<x}|$. Also,
because $P_{<x}\subseteq\minusmax{P}$, we have $P_{<x}=\minusmax{P}$ and
thus 
$P$ is LOT, as desired.
\end{proof}

\begin{definition}  \label{muP}
For any poset $P$, let $\mu(P)$ be defined by
$$
\mu(P)=\max\{ k\ :\ \textrm{ there exists } S \subseteq \max P \textrm{
with } |S|=k\textrm{ such that if } x,y\in S \textrm{ then } P_{<x}=P_{<y}
\}
$$
\end{definition}

Clearly $\mu(P)\le |\max P|$, with equality if and only if $P$ is LOT.  It
turns out that $\mu(P)$ is a useful statistic for induction.  We now have
all 
the necessary tools to prove Theorem~\ref{strongstrom}.

\begin{proofofstrongstrom}
Notice that the claim is trivial for $|P|<l+1$ as all posets on less than
$l+1
$ elements cannot have any $Q_l$-patterns and thus they have the trivial
$l$-
decomposition $P\oplus\emptyset$.

Assume to the contrary that the claim is not true and let $P$ be an $Q_l$-
optimal poset of least cardinality that does not have a LOT
$l$-decomposition 
with $\mu(P)$ maximal over all such choices of $P$ and $|P|\ge l+1$. Let $S$
be the set from Definition~\ref{muP}, $m=|\max P|$ and
$k=|\minusmax{P}|=|P|-
m$.

First, we claim that $P$ is LOT.  If not, then there is some element, say
$x\in (\max P)\setminus S$.  Also let $y\in S$.  If $\packpatabbrev(P;x)\ne
\packpatabbrev(P;y)$, then by Lemma~\ref{stromLOT} either $\packpatabbrev(P^
{x \rightarrow y})>M_{|P|}$ or $\packpatabbrev(P^{y \rightarrow
x})>M_{|P|}$, 
both contradictions, so $\packpatabbrev(P;x) = \packpatabbrev(P;y)$ and
$\pxy$ is $Q_l$-optimal.  Since $\mu(\pxy)>\mu(P)$, by our choice of $P$ we
know that $\pxy$ has an $l$-decomposition $P_1\oplus P_2$.

If $P_2=\emptyset$, then $\packpatabbrev(\pxy)=0$, so by Lemma~\ref
{Mmonotone}, $|P|<l+1$ (because $M_{l+1}=1$), a contradiction to our choice
of $P$.

Hence we may assume that $P_2\ne\emptyset$, so $\pxy$ is LOT.  As the only
element $P$ and $\pxy$ disagree on is $x$, we have that $P-x$ is LOT.
Hence by Lemma~\ref{minusa}, $P$ is also LOT.

Now that we know that $P$ is LOT, we get that $\packpatabbrev(P)
=\packpatabbrev(\minusmax{P})+{m\choose l}k$, so $\minusmax{P}$ is $Q_l$-
optimal.  By induction, $\minusmax{P}$ has an $l$-decomposition
$\minusmax{P}
=P_1\oplus P_2$ and thus $P=P_1\oplus(P_2\oplus\max P)$ is an $l$-
decomposition for $P$.
\end{proofofstrongstrom}

Note that by using the ideas in the last paragraph of  this proof one
may show that if $P=P_1\oplus
P_2$ is an $l$-decomposition for an $Q_l$-optimal poset $P$ then
$|P_1|<l+1$.  Hence because all posets on less than three elements are
layered, all $P_{132}$-optimal posets (and thus $132$-optimal permutations)
are layered.  This observation will be useful in the following proof.
 
\begin{theorem} \label{mn-1}
If $P$ is such that $\packpat(P)=M_{|P|,Q_l}-1$ then there is a poset $Q$
with $|Q|=|P|$ and one of the following:
\ben
\item[(i)] $\packpat(Q)=M_{|P|,Q_l}-1$ and $Q$ is LOT, or
\item[(ii)] $Q$ is $Q_l$-optimal and $|\max Q|=l$, or
\item[(iii)] $l=2$ and $|P|=5.$
\een
\end{theorem}
\begin{proof}
Assume that (i) does not hold and choose $P$ with
$\packpatabbrev(P)=M_{|P|}-1
$ and $\mu(P)$ maximal over all such choices.  Let $n=|P|$, $m=|\max P|$ and
$k=|\minusmax{P}|=n-m$.

We must have $|\packpatabbrev(P;x)-\packpatabbrev(P;y)|\le 1$ for all
$x,y\in\max P$ as otherwise by Lemma~\ref{stromLOT} we would have either
$\packpatabbrev(P^{y\rightarrow x})>M_{n}$ or $\packpatabbrev(P^
{x\rightarrow y})>M_{n}$, a contradiction.  Hence we have
\[\max\{\packpatabbrev(P;x)-\packpatabbrev(P;y) : x,y\in\max
P\}\in\{0,1\}.\]

First we tackle the easier case, where $\max\{\packpatabbrev(P;x)-
\packpatabbrev(P;y) : x,y\in\max P\}=1$.  Pick two maximal elements of $P$,
say $x,y\in\max P$, so that $\packpatabbrev(P;y)-\packpatabbrev(P;x)=1$.  By
Lemma~\ref{stromLOT} we have that $\packpatabbrev(\pxy)=M_{n}$, and thus
by 
Theorem~\ref{strongstrom} we know $\pxy$ has an $l$-decomposition $P_1\oplus
P_2$.  Since $\packpatabbrev(P)=M_{n}-1$, we must have $M_{n}>0$, so we
also have that $n\ge l+1$ and $P_2\ne\emptyset$.  Therefore $\pxy$ and
consequently $P-x$ are LOT.  Hence by Lemma~\ref{minusa}, $P$ is LOT,
a contradiction. 

Now assume $\max\{\packpatabbrev(P;x)-\packpatabbrev(P;y) : x,y\in\max P\}=0
$.  Let $S$ be as in Definition~\ref{muP}, pick $x\in(\max P)\setminus S$
($x$ must exist as $P$ is not LOT) and $y\in S$.  Now $\packpatabbrev(P;x)
=\packpatabbrev(P;y)$ and thus $\packpatabbrev(\pxy)\ge M_{n}-1$ by
Lemma~\ref{stromLOT}.  However if $\packpatabbrev(\pxy)=M_{n}-1$ then we
have contradicted our choice of $P$ as $\mu(\pxy)>\mu(P)$.  Therefore
$\packpatabbrev(\pxy)=M_{n}$ so by Theorem~\ref{strongstrom}, $\pxy$ has
an 
$l$-decomposition $P_1\oplus P_2$.  By the same reasoning as the previous
case, $P_2\ne\emptyset$, so again $\pxy$ and $P-x$ are both LOT.

Although we cannot apply Lemma~\ref{minusa} in this case,
(\ref{for7.9}) still holds for $P$ with $a=0$, so
$$
\packpatabbrev(\pxy)-\packpatabbrev(P)=1={m-2\choose l-2}\left(k-|P_{<x}
|\right).
$$
Therefore we must have ${m-2\choose l-2}=1$.  If $l>2$, this implies that
$m=l$, so (ii) is true with $Q=\pxy$.

If $l=2$ then we must have $k-|P_{<x}|=1$, so there is precisely one
element, say $z\in \minusmax{P}\setminus P_{<x}$.  Since $P-x$ is LOT, $z$
must lie in $\max\minusmax P$.  Let $b=|\max\minusmax P|$.  Then we have
\beq    \label{cP}
\packpatabbrev(P)=\packpatabbrev(\minusmax{P})+\packpatabbrev(P-z;\max P)
+\packpatabbrev(P-x;z,\max P)+\packpatabbrev(P;x,z)
\eeq
Because $P-z$ is LOT, we have that $\packpatabbrev(P-z;\max P)={m\choose
2}(k-1)$, and because $P-x$ is LOT we have that $\packpatabbrev(P-x;z,\max
P)={m-1\choose 2}$.  Notice that because $P^{x\rightarrow y}$ is
$\achain_1\oplus\achain_2$-optimal, by the comment after the proof of
Theorem~\ref{strongstrom}, $P^{x\rightarrow y}$ is layered, and thus
$\minusmax{P}$ is layered.  Since the $\achain_1\oplus\achain_2$-patterns in
$P$ containing both $x$ and $z$ are formed with exactly one element which
lies in $P_{<z}$, $\packpatabbrev(P;x,z)=k-b$.
Finally, $c(\pxy)=c(\bar{P})+{m\choose l}k$.  Now combining all these
$c$-values with equation~\ree{cP} gives
\beq\label{l=2}
\packpatabbrev(\pxy)-\packpatabbrev(P)=1={m\choose 2}k-{m\choose
2}(k-1)-{m-1\choose 2}-k+b,
\eeq
so $k+2=b+m$.  We have by Lemma~\ref{le} (iii) that $m>k$ and $b>k/2$ (this
follows from the fact that $\minusmax{P}$ is layered and
$\achain_1\oplus\achain_2$-optimal), which forces $k\le3$.
This in turn implies $|P|=k+m\le7$.
Now it can be checked by direct computation that for $|P|$ in this
range either the theorem is true vacuously or one of (i) to (iii) holds.
\end{proof}

\begin{theorem}\label{mn-1followup}
If there is an $n$-poset $P$ with $c_{Q_l}(P)=M_{n,Q_l}-1$ then there is an
$n$-poset $Q$ with $c_{Q_l}(Q)=M_{n,Q_l}-1$ and
\ben
\item[(i)] if $l>2$ then $Q$ is layered, or
\item[(ii)] if $l=2$ then $Q=Q_1\oplus Q_2$ where $|Q_1|\le 5$ and $Q_2$ is
layered.
\een
Furthermore, in either case $Q=P_r\oplus\achain_m$ for some
permutation $r\in S_{n-m}$ and integer $m$ which is positive unless
$l=2$ and $n=5$.
\end{theorem}
\begin{proof}
Induct on $n$.
If $n<l+1$, then $M_n=0$, so the theorem is true vacuously.  If
$n=l+1$, then $M_n=1$ and $\packpatabbrev(\achain_{l+1})=0=M_n-1$.  Hence we
may
assume that $n>l+1$.

If case (ii) of Theorem~\ref{mn-1} is true, let $Q$ be the poset guaranteed
there, $k=|\minusmax{Q}|$ and $m=|\max Q|=l$.  Then by
Lemma~\ref{le} (iii), $k<(k+m)/l<2m/m=2$, so $n=k+m\le l+1$,
a case we have already dealt with.

It is routine to check that the poset $P_{15423}$ satisfies case (ii) of
this theorem if case (iii) of Theorem~\ref{mn-1} is true.

Therefore we may assume that case (i) of Theorem~\ref{mn-1} is true, and
thus there is a LOT $n$-poset $Q$ so that $\packpatabbrev(Q)=M_n-1$.
Since $Q$ is LOT,
$Q=\minusmax{Q}\oplus\max Q=\minusmax{Q}\oplus\achain_m$.  As
$\packpatabbrev(Q)=\packpatabbrev(\minusmax{Q})+k{m\choose l}$, we must
have $\packpatabbrev(\minusmax{Q})\ge M_k-1$.  If
$\packpatabbrev(\minusmax{Q})=M_k$, then by Theorem~\ref{layered}, there is
some layered $k$-poset $R$ so that $\packpatabbrev(R)=M_k$, and thus
$R\oplus\achain_m$ is layered, $\packpatabbrev(R\oplus\achain_m)=M_n-1$ and
$R=P_r$ for some $r\in S_k$.  If $\packpatabbrev(\minusmax{Q})=M_k-1$, then
by induction, there is some poset $R$, $|R|=k$, which satisfies
this theorem.  So $R\oplus\achain_m$ is the desired poset.
\end{proof}

\begin{theorem}\label{IZl=2}
For the pattern $q=\permpat$, there are infinitely many IZ integers.
\end{theorem}
\begin{proof}
Assume that the theorem is false.  Since $S_6(M_6-1)=0$ for $l=2$ and
$S_{l+2}(M_{l+2}-1)=0$ for $l\ge3$, there
must be some maximal $k\ge l+2$ so
that $S_k(M_k-1)=0$.  By Theorem~\ref{cont}, there is some $n$ so that
$k_{n-2}=k-1$ and $k_{n-1}=k$.  Also note that since $k_n\ge k_{n-1}=k\ge
l+2$, by Lemma~\ref{le} (iii) we have $n>lk_n\ge l(l+2)>2l+2$, so we
may apply Lemma~\ref{cnidiff} to see that $k_n=k$.

By our choice of $k$, $S_n(M_n-1)\ne 0$, so there is some $p\in S_n$ so that
$\packpatabbrev(p)=M_n-1$.  By Lemma~\ref{permposetequiv},
$\packpatabbrev(P_p)=M_n-1$, and thus Theorem~\ref{mn-1followup}
produces a poset $Q=P_r\oplus\achain_m$ for some $r\in S_{n-m}$ and
integer $m$ which is positive since $n>6$.

Let $\bar{k}=n-m$.  By Theorem~\ref{layered}, there is a layered
$Q_l$-optimal $\bar{k}$-poset $R$, and so we must have
$\packpatabbrev(R\oplus\achain_m)\ge \packpatabbrev(Q)$.  Therefore,
by Lemma~\ref{permposetequiv}, we have
$\packpatabbrev(R\oplus\achain_m)=c_{n,\bar{k}}\ge\packpatabbrev(Q)=M_n-1$,
and thus the inequality in Lemma~\ref{cnidiff} implies that
$\bar{k}=k$.  However, if $\bar{k}=k$ then we
have $\packpatabbrev(r)=M_k-1$, contradicting our choice of $k$.
\end{proof}

Numerical evidence and the contrast between Proposition~\ref{12} and
Theorem~\ref{12...l} amkes us suspect that Theorem~\ref{nonzero} is
not true for $q=\permpat$, $l\ge3$.  In fact, we believe the following
is true.
\begin{conjecture}
The frequency sequence for $q=\permpat$, $l\ge3$ has internal zeros
for all $n\ge l+1$.
\end{conjecture}
It would be interesting to find a proof of this conjecture.  Perhaps a
first step would be to find a simpler proof of Theorem~\ref{IZl=2}.

\noindent{\it Acknowledgement.}  This paper was written, in part, when
Bruce Sagan was resident at the Isaac Newton Institute for
Mathematical Sciences.  He would like to express his appreciation for
the support of the Institute during this period.


\begin{thebibliography}{99}


\bibitem{prec} M.  B\'ona, The number of permutations with exactly $r$
132-subsequences  is P-recursive in the size!, {\it Adv.\ Appl.\ Math.\/}
{\bf 18} (1997), 510--522.


\bibitem{bon:pwo} M. B\'ona, Permutations with one of two
132-subsequences, {\it Discrete Math.\/} {\bf 181} (1998), 267--274.


\bibitem{cw:fsc} T. Chow and J. West, Forbidden subsequences and
Chebyshev polynomials, {\it Discrete Math.\/} {\bf 32} (1980),
125--161.

\bibitem{jr:cfc} M. Jani and R. G. Rieper, Continued fractions and the
Catalan problem, preprint.


\bibitem{kra:prp}  C. Krattenthaler, Permutations with restricted
patterns and Dyck paths, preprint.


\bibitem{man:pca} T. Mansour, Permutations containing and avoiding
certain patterns, preprint.

\bibitem{noo:npc} J. Noonan, The number of permutations containing
exactly one increasing subsequece of length 3, {\it Discrete Math.\/}
{\bf 152} (1996), 307--313.


\bibitem{zeno} J.  Noonan, D.  Zeilberger, The enumeration of
permutations with a prescribed number of ``forbidden'' subsequences,
{\it Adv.\ Appl.\ Math.\/} {\bf 17} (1996), 381--407.


\bibitem{Price} A. Price, ``Packing densities of layered patterns,''
Ph.D. thesis, University of Pennsylvania, Philadelphia, PA, 1997.


\bibitem{rwz:ppc} A. Robertson, H. S. Wilf, and D. Zeilberger,
Permutation patterns and continued fractions, {\it Electronic
J. Combin.\/} {\bf 6} (1999), 6 pages.



\bibitem{ss:rp} R. Simion and F. W. Schmidt, Restricted permutations,
{\it European J. Combin.\/} {\bf 6} (1985), 383--406.


\bibitem{stanley} R. P. Stanley, Log-concave and unimodal
sequences in algebra, combinatorics, and geometry, in
``Graph Theory and Its Applications: East and West,'' {\it Ann.\ NY
Acad.\ Sci.\ } {\bf 576} (1989), 500--535.

\bibitem{sta:ec1} R. P. Stanley, ``Enumerative Combinatorics,
Volume 1,''  Cambridge University Press, Cambridge, 1997.

\bibitem{str:plp} W. Stromquist, Packing layered posets into posets,
manuscript.

\end{thebibliography}
\end{document}